\documentclass[a4paper,12pt,leqno]{article}

\usepackage{amsmath,amssymb,latexsym,amsfonts,amsthm}

\theoremstyle{plain}
\newtheorem{Thm}{Theorem}[section]

\newtheorem{Lem}[Thm]{Lemma}
\newtheorem{Prop}[Thm]{Proposition}

\theoremstyle{definition}
\newtheorem{Def}[Thm]{Definition}
\newtheorem{Rem}[Thm]{Remark}
\newtheorem{Ex}[Thm]{Example}

\newcommand{\bR}{\ensuremath{\mathbb{R}}}

\newcommand{\cE}{\ensuremath{\mathcal{E}}}

\newcommand{\cL}{\ensuremath{\mathcal{L}}}
\newcommand{\cM}{\ensuremath{\mathcal{M}}}

\newcommand{\vn}{\ensuremath{\mbox{{\boldmath $n$}}}}

\newcommand{\vP}{\ensuremath{\mbox{{\boldmath $P$}}}}
\newcommand{\vQ}{\ensuremath{\mbox{{\boldmath $Q$}}}}
\newcommand{\vR}{\ensuremath{\mbox{{\boldmath $R$}}}}

\newcommand{\vU}{\ensuremath{\mbox{{\boldmath $U$}}}}

\newcommand{\eps}{\ensuremath{\varepsilon}}
\newcommand{\e}{{\rm e}}
\renewcommand{\d}{{\rm d}}

\renewcommand{\tilde}{\widetilde}

\newcommand{\cbra}[1]{\left( #1 \right)}
\newcommand{\kbra}[1]{\left\{ #1 \right\}}
\newcommand{\ebra}[1]{\left[ #1 \right]}

\newcommand{\lt}{\ell}

\numberwithin{equation}{section}

\newcounter{No}

\newcounter{Ci}[subsection]
\renewcommand{\theCi}{\arabic{Ci}$^{\circ}$}
\newcommand{\Circ}{\noindent \refstepcounter{Ci} {\bf \theCi).} }

\makeatletter
\renewcommand\section{\@startsection {section}{1}{\z@}%
                                   {-3.5ex \@plus -1ex \@minus -.2ex}%
                                   {2.3ex \@plus.2ex}%
                                   {\normalfont\large\bf}}
\makeatother

\begin{document}

\begin{center}
{\Large \bf Time change approach to generalized excursion measures, 
and its application to limit theorems}
\end{center}
\begin{center}
Patrick J. \textsc{Fitzsimmons} and Kouji \textsc{Yano}\footnotemark \label{foot}
\end{center}
\footnotetext{\ref{foot} E-mail: yano@kurims.kyoto-u.ac.jp}
\bigskip

\begin{center}
{\small Dedicated to Professor M.~Fukushima on the occasion of his seventieth  birthday.}
\end{center}
\bigskip

\begin{abstract}
It is proved that generalized excursion measures can be constructed 
via time change of It\^o's Brownian excursion measure. 
A tightness-like condition on strings is introduced 
to prove a convergence theorem of generalized excursion measures. 
The convergence theorem is applied to obtain a conditional limit theorem, 
a kind of invariance principle where the limit is the Bessel meander. 
\end{abstract}

\section{Introduction}

Stone \cite{MR0158440} has proved various limit theorems for Markov processes 
via time change of Brownian motion. 
The condition assumed is the pointwise convergence of strings. 
Recently Kotani \cite{K1} and Kasahara--Watanabe \cite{KW} have posed 
a tightness-like condition 
to study the scaling limit of the fluctuation of various degenerate limits. 
The second author \cite{Y} has studied this problem 
from the viewpoint of {\em generalized excursion measures} 
through a spectral theoretic approach. 

In the present paper 
we construct generalized excursion measures 
{\em via time change of It\^o's Brownian excursion measure}. 
We introduce a weaker tightness-like condition than Kotani and Kasahara--Watanabe's 
and prove a convergence theorem for generalized excursion measures. 
We apply the convergence theorem to generalize the conditional limit theorem 
obtained by Li--Shiga--Tomisaki \cite{MR2028667}. 

\

\noindent
{\bf (a)} 
For a string $ m $, we construct the generalized excursion measure $ \vn_m $ 
for the $ \cL_m $-diffusion process, where $ \cL_m = \frac{\d}{\d m} \frac{\d}{\d x} $. 
Let $ \vn_{\rm BE} $ denote It\^o's Brownian excursion measure 
and let $ \lt(t,x) $ denote the local time at $ x \in (0,\infty ) $ 
of the excursion path $ (e(t):t \ge 0) $ under $ \vn_{\rm BE} $. 
Define $ A_m(t) = \int_{(0,\infty )} \lt(t,x) \d m(x) $. 
We will show in Lemma \ref{lem: am finite} that 
$ A_m(t) $ is finite or infinite 
$ \vn_{\rm BE} $-almost everywhere, according as $ \int_{0+} x \d m(x) $ is finite or infinite. 
If $ \int_{0+} x \d m(x) < \infty $, then 
we can consider the time-changed excursion path $ e_m(\cdot) = e(A_m^{-1}(\cdot)) $ 
to obtain the desired measure by setting 
$ \vn_m(e \in \cdot) = \vn_{\rm BE}(e_m \in \cdot) $ (Theorem \ref{thm: exc meas}).

We introduce the following condition, 
which we will call the {\em $ \cM_L $-tightness condition}: 
\begin{align}
\lim_{\delta \to 0+} \limsup_{\lambda \to \infty } 
\int_{(0,\delta]} x \log \log (1/x) \d m_{\lambda}(x) = 0 . 
\label{intro cMlog2 tight}
\end{align}
Suppose that $ m_{\lambda}(x) \to m(x) $ for each continuity point $ x $ of $ m $ 
and that the $ \cM_L $-tightness condition is satisfied. 
Then our Theorem \ref{thm: exc conv} asserts that 
the excursion measure for $ \cL_{m_{\lambda}} $ 
converges to that for $ \cL_m $ 
in the sense that 
the time-changed excursion paths under $ \vn_{\rm BE} $ converge uniformly: 
\begin{align}
\lim_{\lambda \to \infty } \sup_{t \ge 0} | e_{m_{\lambda}}(t) - e_m(t) | = 0, 
\qquad \text{$ \vn_{\rm BE} $-a.e.} 
\label{intro unif conv ae}
\end{align}

We will see 
in Proposition \ref{prop: exc conv converse} 
that the following {\em $ \cM $-tightness condition} 
\begin{align}
\lim_{\delta \to 0+} \limsup_{\lambda \to \infty } 
\int_{(0,\delta]} x \d m_{\lambda}(x) = 0 
\label{intro cM tight}
\end{align}
is close to necessary for the convergence \eqref{intro unif conv ae}. 

\

\noindent
{\bf (b)} 
Iglehart \cite{MR0362499} and Bolthausen \cite{MR0415702} 
have proved that 
a suitably rescaled  random walk on $ \bR$, 
conditioned to be positive until a fixed time, 
converges in law to the Brownian meander. 
Li--Shiga--Tomisaki \cite{MR2028667}
have generalized these results to prove that 
a suitably rescaled process for null recurrent generalized diffusion processes 
converges in law to the Bessel meander of dimension $ d \in (0,2) $. 

We generalize the conditional limit theorem of Li--Shiga--Tomisaki 
for positive recurrent diffusion processes. 
Let $ m $ be a string on $ (0,\infty ) $. 
Following Kasahara--Watanabe \cite{KW}, 
we suppose that $ m(\infty )<\infty $ 
and that the difference $ m(\infty ) - m(x) $ is regularly varying as $ x\to\infty $. 
Then a suitably rescaled process $ e^{\lambda,u(\lambda)}(t):= e(\lambda t)/u(\lambda)$, 
conditioned to be positive until time one, namely 
\begin{align}
\cbra{ e^{\lambda,u(\lambda)}(t) : t \in [0,1] } 
\ \text{under} \ 
\vQ^x_m \cbra{ \cdot \mid \zeta(e^{\lambda,u(\lambda)})> 1 } , 
\label{LST}
\end{align}
converges in law to a Bessel meander on $ [0,1] $ of negative dimension. 
See Section \ref{sec: lim thm} for details. 

\

\noindent
{\bf (c)} 
We must remark that 
the importance of certain tightness-like conditions in the class of strings 
was first pointed out by Kotani \cite{K1} and Kasahara--Watanabe \cite{KW}. 
\bigskip

We denote by $ \cM_0 $ 
the class of strings $ m $ for which $0$ is both exit and entrance for $ \cL_m $.
For a family of strings $ \{ m_{\lambda} \} $ in the class $ \cM_0 $ 
such that $ m_{\lambda}(0+) \ge 0 $ for all $\lambda$, 
the pointwise convergence condition 
\begin{align}
\text{$ m_{\lambda}(x) \to m(x) $ \quad as $ \lambda \to \infty$, 
for each continuity point $ x $ of $ m$},
\label{intro ptwise conv}
\end{align}
plays an important role in various limit theorems 
for null recurrent diffusion processes. 
Moreover, this condition is also essential; 
In fact, 
Kasahara \cite{MR0405615} has proved the bi-continuity of Krein's correspondence 
where the class $ \cM_0 $ is equipped with a topology 
induced by the pointwise convergence \eqref{intro ptwise conv}. 
Based on this theory 
Kasahara \cite{MR0405615} and Watanabe \cite{MR1335470}
have established the converse of limit theorems. 

A major breakthrough both in the limit theory and in the spectral theory 
has recently been achieved by Kotani \cite{K1} and Kasahara--Watanabe \cite{KW}. 
Let $ \cM_1 $ denote the class of strings $ m $ 
for which the origin  is of limit circle type for $ \cL_m$. 
They have introduced the following condition 
for strings $ m_{\lambda} \in \cM_1 $: 
\begin{align}
\lim_{\delta \to 0+} \limsup_{\lambda \to \infty } 
\int_0^{\delta} m_{\lambda}(x)^2 \d x = 0 . 
\label{intro cM1 tight}
\end{align}
This condition seems to mean a kind of tightness 
for the family of Radon measures $ \d m_{\lambda}(x) $. 
For this reason 
we call the condition \eqref{intro cM1 tight} the {\em $ \cM_1 $-tightness condition}. 
We say that $ m_{\lambda} \to m $ in $ \cM_1 $ 
if both the pointwise convergence condition \eqref{intro ptwise conv} 
and the $ \cM_1 $-tightness condition \eqref{intro cM1 tight} hold. 
Kotani has generalized 
Kasahara's continuity theorem for Krein's correspondence to the class $ \cM_1 $. 
Kasahara--Watanabe have obtained limit theorems 
for the fluctuation 
of the occupation time and the inverse local time 
of positive recurrent diffusion processes. 

\

\noindent
{\bf (d)} 
Let us recall the usual excursion theory. 
Consider a diffusion process on $[0,\infty)$ for which the origin is a reflecting boundary. 
Note that the origin is necessarily both exit and entrance in Feller's sense. 
Then one can construct the excursion point process 
by using the zero set of a sample path to cut the path into excursions. 
It\^o \cite{MR0402949} has shown that 
the point process is stationary Poisson whose characteristic measure 
is given by a $ \sigma $-finite measure, 
which we call the {\em It\^o excursion measure}. 
Conversely, one can construct a diffusion process 
from an It\^o excursion measure 
by stringing together the excursions of the associated point process. 

The It\^o excursion measures have several descriptions. 
We list the following four formulae 
for the Brownian excursion measure $ \vn_{\rm BE} $. 
See the textbooks \cite{MR1011252}, \cite{MR1725357} and \cite{MR1138461} for details. 

Let $ M $ stand for the maximum value of the path. 
Let $ \vQ_{\rm BM}^x $ for $ x>0 $ stand for the law of the one-dimensional Brownian motion 
(with generator $ \frac{1}{2} \frac{\d ^2}{\d x^2} $) 
starting from $ x $ and absorbed at the origin. 
Let $ \vP_{\rm 3B}^x $ for $ x \ge 0 $ stand for the $ h $-transform of $ \vQ_{\rm BM}^x $ 
with respect to the Brownian scale, 
which is actually the law of the 3-dimensional Bessel process 
(with generator $ \frac{1}{2} \frac{\d ^2}{\d x^2} + \frac{1}{x} \frac{\d}{\d x} $) 
starting from $ x $. 
\\
{\rm (i)} $ \vn_{\rm BE}(M=0)=0 $ and 
for each $ x>0 $ and every bounded continuous functional $ F $ on the excursion space 
whose support is contained in $ \{ M>x \} $ for some $ x>0 $, 
\begin{align}
\vn_{\rm BE} \ebra{F} 
= \lim_{a \to 0+} \frac{1}{a} \vQ_{\rm BM}^a \ebra{F} . 
\label{BE limit}
\end{align}
{\rm (ii)} The strong Markov property: Under $\vn_{\rm BE} $ the excursion process
$ (e(t):t \ge 0) $ is a strong Markov process 
with transition kernel $ \vQ_{\rm BM}^x(e(t) \in \d y) $ 
and entrance law $ \frac{1}{x} \vP_{\rm 3B}^0(e(t) \in \d x) $. 
In particular, for each positive stopping time $ \tau $ 
and every measurable set $ \Gamma $, 
\begin{align}
\vn_{\rm BE} (e(\tau+\cdot) \in \Gamma ) 
= \int_{(0,\infty )} \frac{1}{x} \vP_{\rm 3B}^0(e(\tau) \in \d x) 
\vQ_{\rm BM}^x(\Gamma ) . 
\label{BE sMarkov}
\end{align}
\\
{\rm (iii)} The maximum decomposition due to Williams \cite{MR0350881}: 
For any measurable set $ \Gamma $, 
\begin{align}
\vn_{\rm BE}(\Gamma ) = \int_0^{\infty } \vR_{\rm 3B}^x(\Gamma ) \frac{\d x}{x^2}. 
\label{BE Williams}
\end{align}
Here $ \vR_{\rm 3B}^x $ stands for the law of the path 
defined by piecing together two independent $ \vP_{\rm 3B}^x $-processes (the second run backward in time)
until they first hit the point $ x $. 
The formula \eqref{BE Williams} means that 
the law of the maximum $ M $ is given as $ \vn_{\rm BE}(M \in \d x) = \frac{\d x}{x^2} $ 
and the conditional law of $ \vn_{\rm BE} $ given $ M=x $ is $ \vR_{\rm 3B}^x $. 
\\
{\rm (iv)} The lifetime decomposition: 
For any measurable set $ \Gamma $, 
\begin{align}
\vn_{\rm BE}(\Gamma ) = \int_0^{\infty } \vP_{\rm 3B}^{0}(\Gamma |e(t)=0) p_{\rm 3B}(t,0,0) \d t. 
\label{BE lifetime}
\end{align}
Here $ p_{\rm 3B}(t,x,y) $ stands for the transition probability density of $ \vP^x_{3B} $ 
with respect to its speed measure. 
The formula \eqref{BE lifetime} means that 
the law of the lifetime $ \zeta $ is given 
as $ \vn_{\rm BE}(\zeta \in \d t) = p_{\rm 3B}(t,0,0) \d t $ 
and the conditional law of $ \vn_{\rm BE} $ given $ \zeta=t $ 
is $ \vP_{\rm 3B}^{0}(\cdot|e(t)=0) $. 

\

\noindent
{\bf (e)} 
Let $ m $ be a string on $ (0,\infty ) $ 
such that the origin is exit and non-entrance for $ \cL_m $.  
Then there is no  reflecting $ \cL_m $-diffusion process 
and therefore the usual excursion theory is not available. 
Nevertheless 
one can study a $ \sigma $-finite measure on the excursion space 
such that the description formulae listed above still holds 
where we replace $ \vQ_{\rm BM}^x $ by the absorbing $ \cL_m $-diffusion process 
and $ \vP_{\rm 3B}^x $ by its $ h $-transform, etc. 
We call such a measure the {\em generalized excursion measure for the $ \cL_m $-diffusion process}. 

Pitman--Yor \cite{MR656509} have introduced such generalized excursion measures 
and established the description formulae {\rm (i)}, {\rm (ii)} and {\rm (iii)}. 
In \cite{MR1439532} they established the formula {\rm (iv)} 
for the Bessel processes of dimension $ -\infty < d < 2 $. 
They used the generalized excursion measures for Bessel processes 
to obtain several remarkable path decomposition formulae 
for Bessel processes and Bessel bridges. 

The use of time change, in the context of It\^o's Brownian excursion measure, 
can be found in Biane--Yor \cite{MR886959}, Section 3 (See also \cite{MR799465}). 
Their motivation was to compute the joint law of $ (H(\eta(t)),\eta(t)) $ 
where $ H(t) $ denotes Cauchy's principal value of Brownian motion 
and $ \eta(t) $ denotes the inverse local time at the origin. 

The second author \cite{Y} has established 
the lifetime decomposition formula {\rm (iv)} 
for generalized excursion measures 
assuming that 
the Laplace transform of the absorbing spectral measure, 
which is the counterpart of $ p_{\rm 3B}(t,0,0) $ in \eqref{BE lifetime}, 
is finite. 
He also studied the relationship between the absorbing spectral measure 
and the (reflecting) spectral measure corresponding to {\em dual strings} 
to study Kasahara--Watanabe's limit theorem from this viewpoint. 
We do not go into the lifetime decomposition formula 
in the present paper. 

\

The paper is organized as follows. 
In Section \ref{sec: results} we state our results. 
We will state the construction theorem of generalized excursion measures 
in Section \ref{sec: construction} 
and the convergence theorem of generalized excursion measures 
in Section \ref{sec: continuity exc}. 
Applications to limit theorems will be stated in Section \ref{sec: lim thm}. 
Sections \ref{sec: prf construction}, \ref{sec: prf continuity} and \ref{sec: prf lim thm} 
are devoted to the proofs of the results 
given in Sections \ref{sec: construction}, \ref{sec: continuity exc} and \ref{sec: lim thm}, 
respectively.


\vfill\eject

\section{Results} \label{sec: results}

Before stating our results 
we prepare some notation. 

Let $ E $ be the space of continuous paths $ e:[0,\infty ) \to [0,\infty ) $ 
such that if $ e(t_0)=0 $ for some $ t_0>0 $ then $ e(t)=0 $ for all $ t > t_0 $. 
We set 
\begin{align}
\zeta(e) = \inf \kbra{ t>0 : e(t)=0 } 
\end{align}
where $ \inf \emptyset = \infty $. 
We call $ \zeta(e) \in [0,\infty ] $ the {\em lifetime} of the path $ e $. 
Hence, if $ \zeta(e) \in (0,\infty ) $, then 
$ e(t)>0 $ for $ 0<t<\zeta(e) $ and $ e(t)=0 $ for $ t \ge \zeta(e) $. 
We regard $ E $ as a complete separable metric space 
equipped with the compact uniform topology. 
Let $ \cE $ denote its Borel $ \sigma $-field. 
Then all the measures 
$ \vn_{\rm BE} $, $ \vP_{\rm 3B}^x $ for $ x \ge 0 $, 
$ \vR^x_{\rm 3B} $ and $ \vQ^x_{\rm BM} $ for $ x>0 $ 
may be considered to be defined on $ E $. 

For $ x \in (0,\infty ) $, 
denote by $ \tau_x $ the first passage time to $ x $. 
Denote $ M(e) = \max_{t \ge 0} e(t) $. 

We fix versions of the local time at $ x \in (0,\infty ) $ 
for almost every path $ e $ under 
$ \vn_{\rm BE} $, $ \vP_{\rm 3B}^x $ for $ x \ge 0 $, $ \vR^x_{\rm 3B} $ 
and $ \vQ^x_{\rm BM} $ for $ x>0 $. 
We will denote them by the common symbol $ \lt(t,x) $. 
Thus, 
$ \lt(t,x) $ is jointly continuous on $ (0,\infty ) \times (0,\infty ) $ 
and the equality 
\begin{align}
\int_0^t f(e(s)) \d s = 
2 \int_{(0,\infty )} f(x) \lt(t,x) \d x 
\label{local time}
\end{align}
holds for every bounded Borel function $ f $ on $ (0,\infty ) $ 
for almost every path $ e $ with respect to the measure $ \vn_{\rm BE} $, and so on. 
For instance, 
we know by the the maximum decomposition formula \eqref{BE Williams} that 
the process $ (\lt(t,\cdot):t \le \tau_x) $ under $ \vR^x_{\rm 3B} $ for $ x>0 $ 
has the same law as $ (\lt(t,\cdot):t \le \tau_x) $ under $ \vP_{\rm 3B}^0 $.

\subsection*{Classes of strings} 

A  {\em string} $m$ on $(0,\infty )$ is a function $ m:(0,\infty ) \to (-\infty , \infty ) $ 
which is strictly-increasing and right-continuous. 
\begin{Rem}
In the context of generalized diffusion processes, 
strings are only assumed to be non-decreasing and right-continuous. 
\end{Rem}

We consider the following four classes of strings: 
\begin{align}
\cM_0 =& \kbra{ m:\text{string, $ m(0+) $ is finite} }, 
\\
\cM_1 =& \kbra{ m:\text{string}, \ \int_{0+} m(x)^2 \d x < \infty }, 
\\
\cM_L =& \kbra{ m:\text{string}, \ \int_{0+} x \log \log (1/x) \,\d m(x) < \infty } 
\intertext{and} 
\cM =& \kbra{ m:\text{string}, \ \int_{0+} x \,\d m(x) < \infty } . 
\end{align}
Then the following relation holds: 
\begin{align}
\cM_0 \subset \cM_1 \subset \cM_L \subset \cM . 
\end{align}
The relation $ \cM_1 \subset \cM_L $ 
follows from the fact that 
there exists $ C $ such that 
\begin{align}
\int_{(0,\delta]} x \log \log (1/x) \d m(x) 
\le C \cbra{ \int_0^{\delta} m(x)^2 \d x }^{1/2} 
\qquad m \in \cM_1 , \ \delta < 1/2 . 
\label{cM1 subset cMlog2}
\end{align}
In fact, we integrate the LHS by parts to obtain 
\begin{align}
\int_{(0,\delta]} x \log \log (1/x) \d m(x) 
\le C' \int_0^{\delta} |m(x)| \log \log(1/x) \d x 
\qquad m \in \cM_L , \ \delta < 1/2 
\label{dominate}
\end{align}
for some constant $ C' $. 

\begin{Rem}
{\rm (i)} 
The origin is both exit and entrance for $ \cL_m $ 
if and only if $ m \in \cM_0 $. \\
{\rm (ii)} 
The origin is of limit circle type for $ \cL_m $ in Weyl's sense 
if and only if $ m \in \cM_1 $. \\
{\rm (iii)} 
The origin is exit for $ \cL_m $ 
if and only if $ m \in \cM $. 
\end{Rem}

\begin{Ex}
For $ \alpha \in (0,\infty ) $, let 
\begin{align}
m^{(\alpha )}(x) = 
\begin{cases}
(1-\alpha )^{-1} x^{\frac{1}{\alpha} -1} \qquad & \text{if} \ \alpha \in (0,1) , \\
\log x & \text{if} \ \alpha = 1 , \\
- (\alpha -1)^{-1} x^{\frac{1}{\alpha} -1} \qquad & \text{if} \ \alpha \in (1,\infty ) . 
\end{cases}
\label{m alpha}
\end{align}
Note that $ \d m^{(\alpha )}(x) = \alpha^{-1} x^{\frac{1}{\alpha} -2 } \d x $ 
for all $ \alpha \in (0,\infty ) $. 
Then
\itemindent=20pt
\item{(i)} $ m^{(\alpha )} \in \cM_0 $ if and only if $ \alpha \in (0,1) $;

\item{(ii)} $ m^{(\alpha )} \in \cM_1 $ if and only if $ \alpha \in (0,2) $;

\item{(iii)} $ m^{(\alpha )} \in \cM_L $ for all $ \alpha \in (0,\infty ) $. 
\end{Ex}

\subsection{Construction of generalized excursion measures} \label{sec: construction}

For a string $ m $ on $ (0,\infty ) $, we define 
\begin{align}
A_m(t) = \int_{(0,\infty )} \lt(t,x) \d m(x) 
\qquad \text{for $ t \ge 0 $}. 
\label{am}
\end{align}

\begin{Lem} \label{lem: am finite}
Let $ m $ be a string on $ (0,\infty ) $. Then the following dichotomy holds: \\
{\rm (i)} If $ m \in \cM $, then $ A_m(t) < \infty $ for all $ t \ge 0 $, 
$ \vn_{\rm BE} $-a.e. \\
{\rm (ii)} If $ m \notin \cM $, then $ A_m(t) = \infty $ for all $ t \ge 0 $, 
$ \vn_{\rm BE} $-a.e. 
\end{Lem}

In what follows we assume that $ m \in \cM $. 
It is easy to see that the function $ A_m $ is 
continuous and strictly increasing on $ (0,\zeta] $ 
and is constant on $ [\zeta,\infty ) $. 
Define $ A_m^{-1}(t) $ for $ t < A_m(\zeta) $ by the inverse function of $ A_m $ 
and set $ A_m^{-1}(t) = \zeta $ for $ t \ge A_m(\zeta) $. 
Then $ A_m^{-1}(A_m(t)) = t \wedge \zeta $ for all $ t \ge 0 $. 
We define the time-changed process on the space $ E $ by 
\begin{align}
e_m(t) = e(A_m^{-1}(t)) \qquad t \ge 0 
\label{em}
\end{align}
and define a $ \sigma $-finite measure on $ E $ by 
\begin{align}
\vn_m(\cdot) = \vn_{\rm BE}(e_m \in \cdot) . 
\end{align}

For $ x>0 $ let  $ \vQ_m^x $ denote the law of the $ \cL_m $-diffusion process 
starting from $ x $ and absorbed at the origin. 
Let $ \vP_{m}^x $ for $ x \ge 0 $ denote the $ h $-transform of $ \vQ_m^x $ 
with respect to $h(x) =x$ (the Brownian scale function). 
This is actually the law of the  diffusion process (starting from $ x $) with speed measure
$ x^2 \d m(x) $ and scale function $ -1/x $. 
Let $ \vR_{m}^x $ for $ x>0 $ denote the law of the path 
defined by piecing together two independent $ \vP_{m}^x $-processes 
until they first hit the point $ x $ (the second one being run backward in time). 

Then we obtain the following description formulae for $ \vn_m $. 

\begin{Thm} \label{thm: exc meas}
Let $ m \in \cM $. Then the following hold: \\
{\rm (i)}$ ' $ $ \vn_m(M=0)=0 $ and 
for every bounded continuous functional $ F $ on $ E $ 
whose support is contained in $ \{ M>x \} $ for some $ x>0 $, 
\begin{align}
\vn_m \ebra{F} 
= \lim_{a \to 0+} \frac{1}{a} \vQ_m^a \ebra{F} . 
\label{GE limit}
\end{align}
{\rm (ii)}$ ' $ For each positive stopping time $ \tau $ and every measurable set $ \Gamma $, 
\begin{align}
\vn_m (e(\tau+\cdot) \in \Gamma ) 
= \int_{(0,\infty )} \frac{1}{x} \vP_{m}^0(e(\tau) \in \d x) 
\vQ_m^x(\Gamma ) . 
\label{GE sMarkov}
\end{align}
{\rm (iii)}$ ' $ For any measurable set $ \Gamma $, 
\begin{align}
\vn_m(\Gamma ) = \int_0^{\infty } \vR_{m}^x(\Gamma ) \frac{\d x}{x^2}. 
\label{GE Williams}
\end{align}
\end{Thm}

For later use we note that
if we put $ \tau=\tau_x $ (for some fixed  $ x>0 $)
then the formula {\rm (ii)}$ ' $ becomes 
\begin{align}
\vQ_m^x(\Gamma ) 
= x \vn_m (e(\tau_x+\cdot) \in \Gamma ) . 
\label{Qm formula}
\end{align}

\subsection{Convergence theorem of generalized excursion measures} \label{sec: continuity exc}

Let $ \{ m_{\lambda} \} $ be a family of strings on $ (0,\infty ) $. 
As we have mentioned in the introduction, 
the pointwise convergence condition 
\begin{align}
\text{$ m_{\lambda}(x) \to m(x) $ as $ \lambda \to \infty $ 
for  all continuity points  $ x $ of $ m $} , 
\label{ptwise conv}
\end{align}
is inadequate for the studies of various limit theorems and of the spectral theory. 

In the class $ \cM_0 $, it is usually assumed that $ m_{\lambda}(0+) \ge 0 $. 
We may regard it as the tightness condition of the class $ \cM_0 $. 

For the class $ \cM_1 $, 
Kotani and Kasahara-Watanabe have introduced the condition 
\begin{align}
\lim_{\delta \to 0+} \limsup_{\lambda \to \infty } 
\int_0^{\delta} m_{\lambda}(x)^2 \d x = 0 , 
\label{cM1 tight}
\end{align}
which we call the {\em $ \cM_1 $-tightness condition}. 

For the classes $ \cM_L $ and $ \cM $, 
we consider the following tightness-like conditions. 

\begin{Def} \label{def: cM-conv}
Let $ m_{\lambda}, m \in \cM $. \\
{\rm (i)} The condition 
\begin{align}
\lim_{\delta \to 0+} \limsup_{\lambda \to \infty } 
\int_{(0,\delta]} x \log \log (1/x) \d m_{\lambda}(x) = 0 
\label{cMlog2 tight}
\end{align}
is called the {\em $ \cM_L $-tightness condition}. 
\\
{\rm (ii)} The condition 
\begin{align}
\lim_{\delta \to 0+} \limsup_{\lambda \to \infty } 
\int_{(0,\delta]} x \d m_{\lambda}(x) = 0 
\label{cM tight}
\end{align}
is called the {\em $ \cM $-tightness condition}. 
\end{Def}

\begin{Def}
Let $ m_{\lambda}, m \in \cM $. 
We say that $ m_{\lambda} \to m $ in $ \cM_L $ (resp. $ \cM $) 
if both the pointwise convergence condition \eqref{ptwise conv} 
and the $ \cM_L $- (resp. $ \cM $-) tightness condition hold. 
\end{Def}

\begin{Rem}
It is immediate by \eqref{cM1 subset cMlog2} 
that $ \cM_1 $-tightness implies $ \cM_L $-tightness. 
It is obvious by definition 
that $ \cM_L $-tightness implies $ \cM $-tightness. 
\end{Rem}

The following theorem asserts that 
$ \cM_L $-convergence implies 
pathwise uniform convergence of time-changed excursion processes. 
\begin{Thm}[Convergence theorem of generalized excursion measures] \label{thm: exc conv}
Assume that $ m_{\lambda} \to m $ in $ \cM_L $ as $ \lambda \to \infty $. 
Then 
\begin{align}
\lim_{\lambda \to \infty } \sup_{t \ge 0} 
\left| A_{m_{\lambda}}(t) - A_m(t) \right| = 0 
\qquad \text{$ \vn_{\rm BE} $-a.e.} 
\label{am conv}
\end{align}
and 
\begin{align}
\lim_{\lambda \to \infty } \sup_{t \ge 0} | e_{m_{\lambda}}(t) - e_m(t) | = 0 
\qquad \text{$ \vn_{\rm BE} $-a.e.} 
\label{unif conv ae}
\end{align}
\end{Thm}

The following proposition asserts that 
the $ \cM $-tightness condition is necessary for the convergence \eqref{am conv}. 
\begin{Prop} \label{prop: exc conv converse}
Let $ m_{\lambda},m \in \cM $. 
Assume that the convergence \eqref{am conv} holds. 
Then the $ \cM $-tightness condition \eqref{cM tight} holds. 
\end{Prop}

\subsection{Conditional limit theorem} \label{sec: lim thm}

For two functions $f$ and $g$ defined for all large reals, we write  
$ f(x) \sim g(x) $ as $ x \to \infty $ to mean that $ \lim_{x \to \infty } f(x)/g(x)=1 $.
A function $ K $ defined for large (real) $ x $ 
is \emph{slowly varying} as $ x \to \infty $ 
provided $ K(\lambda x) \sim K(\lambda) $ as $ \lambda \to \infty $ for all sufficiently large $ x$. 

For $ e \in E $, $ \lambda_1>0 $ and $ \lambda_2>0 $, we use the notation 
\begin{align}
e^{\lambda_1,\lambda_2}(t) = \frac{1}{\lambda_2} e(\lambda_1 t) 
\qquad \text{for $ t \ge 0 $}. 
\label{rescaled}
\end{align}

First of all, we restate the conditional limit theorem of Li-Shiga-Tomisaki,  using our terminology. 
Consider the following assumption: \medskip
\\
{\bf (A1)} $ \alpha \in (0,1) $, 
$ m \in \cM_0 $ 
and 
\begin{align}
m(x) \sim (1-\alpha )^{-1} x^{\frac{1}{\alpha }-1} K(x) 
\qquad \text{as $ x \to \infty  $,} 
\end{align}
for some  function $ K$ slowly varying at infinity. 

\begin{Rem}
Note that, if the assumption {\bf (A1)} is satisfied, 
then the corresponding $ \cL_m $-diffusion process is null recurrent; in fact, $ m(\infty )=\infty $. 
\end{Rem}

Let $ u(\lambda) $ be an asymptotic inverse of $ \lambda^{\frac{1}{\alpha }} K(\lambda) $. 

\begin{Thm}[Li--Shiga--Tomisaki {\cite[Theorem 3.1]{MR2028667}}] \label{thm: LST}
Suppose that {\bf (A1)} is satisfied. 
Then, for any $ x>0 $, the distribution on $ E $ of 
the rescaled process $ (e^{\lambda,u(\lambda)}(t):t \in (0,1]) $ 
under $ \vQ_m^x(\cdot \mid \zeta(e^{\lambda,u(\lambda)})>1) $ 
converges as $ \lambda \to \infty $ to 
the process $ (e(t):t \in (0,1]) $ 
under $ \vn_{m^{(\alpha )}}(\cdot \mid \zeta>1) $ 
where $ m^{(\alpha )} \in \cM_0 $ has been introduced in \eqref{m alpha}. 
\end{Thm}

\begin{Rem}
For $ \alpha \in (0,\infty ) $, 
the process $ (e(t):t \in [0,1]) $ 
under $ \vn_{m^{(\alpha )}}(\cdot \mid \zeta>1) $ 
is called the Bessel meander of dimension $ d=2-2\alpha $. 
Theorem \ref{thm: LST} says that 
the limit process is a Bessel meander of positive dimension. 
\end{Rem}

\begin{Rem}
More precisely, Theorem 3.1 of \cite{MR2028667} 
covers the case where a string is non-decreasing but is not strictly increasing 
so that the $ \cL_m $-diffusion process is a generalized diffusion process.
\end{Rem}

We introduce the following assumptions on a string $ m \in \cM $: 
\bigskip
\\
{\bf (A2)} $ \alpha =1 $, 
$ m \in \cM_0 $ 
and 
\begin{align}
m(\lambda x) - m(\lambda) \sim (\log x) K(\lambda) 
\qquad \text{as $ \lambda \to \infty $, for all $ x>0 $},
\end{align}
for some function $ K$ slowly varying  at infinity. 
\bigskip
\\
{\bf (A3)} $ \alpha \in (1,\infty ) $. For small $x$, the string $ m $ satisfies 
\begin{align}
\lim_{x \to 0+} \frac{m(\infty ) - m(x)}{x^{\frac{1}{\alpha }-1} } < \infty . 
\label{A3 prime at 0}
\end{align}
For large $ x $ it satisfies 
$ m(\infty )< \infty $ and 
\begin{align}
m(\infty ) - m(x) \sim (\alpha -1)^{-1} x^{\frac{1}{\alpha }-1} K(x) 
\qquad \text{as $ x \to \infty $,} 
\label{A3 prime at infty}
\end{align}
for some  function $ K$ slowly varying at infinity.
\bigskip

\begin{Rem}
Note that the assumption {\bf (A3)} implies that 
the corresponding $ \cL_m $-diffusion process is positive recurrent; 
In fact, $ m(\infty )<\infty $. 
\end{Rem}

Now let us state our generalization of Theorem \ref{thm: LST}. 
Let $ u $ be an asymptotic inverse of $\lambda\mapsto \lambda^{\frac{1}{\alpha }} K(\lambda) $. 

\begin{Thm}[Conditional limit theorem] \label{thm: lim thm}
Suppose that any one of {\bf (A1)}, {\bf (A2)} and {\bf (A3)} is satisfied. 
Then, for each $ x>0 $ and every bounded continuous functional $ F $ on $ E $, 
\begin{align}
\vQ_m^x \ebra{ F(e^{\lambda,u(\lambda)}) \mid \zeta(e^{\lambda,u(\lambda)})>1 } 
\to \vn_{m^{(\alpha )}} \ebra{ F \mid \zeta>1 } 
\qquad \text{as $ \lambda \to \infty $} 
\end{align}
where $ m^{(\alpha )} \in \cM $ has been introduced in \eqref{m alpha} 
and where $ e^{\lambda,u(\lambda)}(t) = e(\lambda t)/u(\lambda) $ 
as is defined in \eqref{rescaled}. 
\end{Thm}

Theorem \ref{thm: lim thm} says that 
in positive recurrent cases 
the limit process is the Bessel meander of the appropriate negative dimension.

The keys to the proof of Theorem \ref{thm: lim thm} are the following three lemmas. 

\begin{Lem} \label{lem: lim thm}
Suppose that any one of {\bf (A1)}, {\bf (A2)} and {\bf (A3)} is satisfied. 
For $ \lambda>0 $, set 
\begin{align}
m_{\lambda}(x) = 
\begin{cases}
m(\lambda x) / \{ \lambda^{\frac{1}{\alpha }-1} K(\lambda) \} 
\qquad & \text{if $ \alpha \in (0,1) $,} 
\\
\kbra{m(\lambda x)-m(\lambda)} / K(\lambda) 
\qquad & \text{if $ \alpha =1 $,} 
\\
\kbra{m(\lambda x)-m(\infty )} / \{ \lambda^{\frac{1}{\alpha }-1} K(\lambda) \} 
\qquad & \text{if $ \alpha \in (1,\infty ) $.} 
\end{cases}
\label{m lambda}
\end{align}
so that $ \d m_{\lambda}(x) = \d m(\lambda x) / \{ \lambda^{\frac{1}{\alpha }-1} K(\lambda) \} $ 
in all cases. 
Then there exist constants $ C $, $ \lambda_0 $ and $ 0<\eps<1 $ such that 
\begin{align}
|m_{\lambda}(x)| \le C x^{\eps -1} 
\qquad \text{for all $ x \in (0,1] $ and for all $ \lambda>\lambda_0 $}. 
\label{global bound}
\end{align}
In particular, 
$ m_{\lambda} \to m^{(\alpha )} $ in $ \cM_L $ as $ \lambda \to \infty $. 
\end{Lem}

\begin{Lem} \label{lem: lim thm2}
Let $ m \in \cM $ 
and suppose that $ \d m_{\lambda}(x) = \lambda \d m(\lambda x) / v(\lambda) $
Then, for each $ x>0 $ and every bounded continuous functional $ F $ on $ E $, 
\begin{align}
\begin{split}
& \vQ_m^x \ebra{ F(e^{v(\lambda),\lambda}) \mid \zeta(e^{v(\lambda),\lambda}) > 1 } 
\\
&\phantom{bbbbb}=\vn_{\rm BE} \ebra{ 
F \cbra{ e(A_{m_{\lambda}}^{-1}(t+A_{m_{\lambda}}(\tau_{x/\lambda}))) : t \ge 0 } 
\mid A_{m_{\lambda}}(\zeta) - A_{m_{\lambda}}(\tau_{x/\lambda}) > 1 } . 
\end{split}
\label{lem: lim thm2 eq1}
\end{align}
Here $ e^{v(\lambda),\lambda}(t) = \lambda^{-1} e(v(\lambda)t) $ 
as is defined in \eqref{rescaled}. 
\end{Lem}

\begin{Lem} \label{lem: lim thm3}
Suppose that $ m_{\lambda} \to m $ in $ \cM_L $ as $ \lambda \to \infty $. 
Suppose, in addition, 
that $ \vn_{\rm BE}(A_m(\zeta)=1)=0 $ 
and that the inequality \eqref{global bound} is satisfied 
for some constants $ C $, $ \lambda_0 $ and $ 0<\eps<1 $. 
Then, for each $ x>0 $ and every bounded continuous functional $ F $ on $ E $, 
\begin{align}
\begin{split}
& \vn_{\rm BE} \ebra{ F(e_{m_{\lambda}}) ; 
A_{m_{\lambda}}(\zeta) - A_{m_{\lambda}}(\tau_{x/\lambda}) > 1 } 
\\
& \phantom{bbbbbbbbb}\to 
\vn_{\rm BE} \ebra{ F(e_m) ; A_m(\zeta) > 1 } 
\qquad \text{as $ \lambda \to \infty $}. 
\end{split}
\label{lem: lim thm3 eq1}
\end{align}
\end{Lem}

We point out that $ \vn_{\rm BE} $ is not a finite measure 
and hence that the bounded convergence theorem does not apply.

\section{Proof of the construction theorem} \label{sec: prf construction}

We begin with the proof of  Lemma \ref{lem: am finite}. 

\begin{proof}[Proof of Lemma \ref{lem: am finite}]
We want to show that 
$ \vn_{\rm BE}(A_m(t)=\infty, \ \exists t \ge 0) = 0 $ if $ m \in \cM $ 
and 
$ \vn_{\rm BE}(A_m(t)<\infty, \ \exists t \ge 0) = 0 $ if $ m \notin \cM $. 

Since $ A_m(t) $ is increasing on $ [0,\zeta] $ and is constant on $ [\zeta,\infty ) $, 
we have 
\begin{align}
\vn_{\rm BE}(A_m(t)=\infty, \ \exists t \ge 0) 
\le \vn_{\rm BE}(A_m(\zeta)=\infty ) . 
\label{am finite eq0}
\end{align}
Since the local time is an additive functional, 
we have $ A_m(\zeta) = A_m(\tau_M) + A_m^{\vee}(\tau_M^{\vee}) $ 
where $ A_m^{\vee} $ (resp. $ \tau_M^{\vee} $) 
is the counterpart of $ A_m $ (resp. $ \tau_M $) 
for the time reversal path. 
Note that $ \{ A_m(\zeta)<\infty \} 
= \{ A_m(\tau_M)<\infty \} \cap \{ A_m^{\vee}(\tau_M^{\vee}) < \infty \} $. 
By definition of $ \vR^a_{\rm 3B} $, 
the random variables $ A_m(\tau_M) $ and $ A_m^{\vee}(\tau_M^{\vee}) $ 
under $ \vR^a_{\rm 3B} $ for $ a>0 $ 
are independent and both have the same distribution as 
$ A_m(\tau_a) $ under $ \vP^0_{\rm 3B} $. 
Hence we obtain 
\begin{align}
\vn_{\rm BE}(A_m(\zeta)=\infty ) 
= \int_0^{\infty } \kbra{ 1 - \vP_{\rm 3B}^0(A_m(\tau_a)<\infty )^2 } \frac{\d a}{a^2} . 
\label{am finite eq1}
\end{align}

Suppose that $ A_m(t)< \infty $ for some $ t \ge 0 $. 
If $ t \ge \tau_M $, then $ A_m(\tau_M) \le A_m(t) < \infty $. 
If $ t<\tau_M $, then $ \inf_{s \in [t,\tau_M]} e(s) > 0 $ and hence 
$ A_m(\tau_M) = A_m(\tau_M)-A_m(t)+A_m(t) < \infty $. 
Thus 
\begin{align}
\vn_{\rm BE}(A_m(t)<\infty, \ \exists t \ge 0) 
& \le\vn_{\rm BE}(A_m(\tau_M)<\infty ) 
\\
& =\int_0^{\infty } \vP_{\rm 3B}^0(A_m(\tau_a)<\infty ) \frac{\d a}{a^2} . 
\label{am finite eq4}
\end{align}

The following 0-1 law is well-known (See, e.g. \cite{MR1071537}, Corollary 1): 
For any $ a>0 $, 
\begin{align}
\vP_{\rm 3B}^0(A_m(\tau_a)<\infty ) = 
\begin{cases}
1 \qquad & \text{if $ m \in \cM $}, 
\\
0 \qquad & \text{if $ m \notin \cM $}. 
\end{cases}
\label{am finite eq3}
\end{align}
Combining \eqref{am finite eq0}-\eqref{am finite eq4} 
with \eqref{am finite eq3}, 
we complete the proof. 
\end{proof}

\begin{Rem}
The 0-1 law \eqref{am finite eq3} is equivalent to 
\begin{align}
\vP_{\rm 3B}^0 \cbra{ \int_{(0,a]} \ell(\infty ,x) \d m(x) <\infty } = 
\begin{cases}
1 \qquad & \text{if $ m \in \cM $}, 
\\
0 \qquad & \text{if $ m \notin \cM $}. 
\end{cases}
\end{align}
This is closely related to Jeulin's lemma (\cite{MR604176}, Lemme (3.22)), 
a useful version of which can be found in Pitman--Yor \cite{MR868512}, Lemma 2. 
\end{Rem}

It is known that, if $ m \in \cM $, then 
$ \vP_{\rm 3B}^x(A_m(t)<\infty , \ \forall t \ge 0) =1 $ for $ x \ge 0 $ and 
$ \vQ_{\rm BM}^x(A_m(t)<\infty , \ \forall t \ge 0) =1 $ for $ x>0 $. 
By a standard time-change argument (as in Chapter 5 of It\^o--McKean's book \cite{MR0345224}), 
we obtain the following. 

\begin{Prop} \label{prop: time-change}
{\rm (i)} 
The time-changed process $ e_m(t)=e(A_m^{-1}(t)) $ under $ \vP_{\rm 3B}^x $ for $ x \ge 0 $ 
is a diffusion process whose law is $ \vP_{m}^x $: 
\begin{align}
\vP_{\rm 3B}^x (e_m \in \cdot) = \vP_{m}^x (\cdot) 
\qquad \text{for all $ x \ge 0 $.} 
\label{time change 1}
\end{align}
{\rm (ii)} 
The time-changed process $ e_m(t)=e(A_m^{-1}(t)) $ under $ \vQ_{\rm BM}^x $ for $ x > 0 $ 
is a diffusion process whose law is $ \vQ_{m}^x $: 
\begin{align}
\vQ_{\rm BM}^x (e_m \in \cdot) = \vQ_{m}^x (\cdot) 
\qquad \text{for all $ x>0 $.} 
\label{time change 2}
\end{align}
\end{Prop}
\noindent The proof is straightforward, so we omit it.

\begin{proof}[Proof of Theorem \ref{thm: exc meas}]
Since $ \tau=A_m^{-1}(t) $ is a positive stopping time, 
we can apply the description formula \eqref{BE sMarkov} to obtain 
\begin{align}
\vn_{\rm BE}(e(A_m^{-1}(t)+\cdot) \in \cdot) 
= \int_{(0,\infty )} \frac{1}{x} \vP^0_{\rm 3B}(A_m^{-1}(t) \in \d x) \vQ_{\rm BM}^x(\cdot) . 
\label{prf exc meas eq1}
\end{align}
Consider the shifted path $ e^+(s) = e(A_m^{-1}(t)+s) $, $ s \ge 0 $. 
We denote the counterpart of $ A_m $ for the shifted path $ e^+ $ by $ A_m^+ $. 
Then we have $ A_m^+(\cdot) = A_m(A_m^{-1}(t)+\cdot) - A_m(A_m^{-1}(t)) $ 
and then we have $ (A_m^+)^{-1}(s) = A_m^{-1}(t+s)-A_m^{-1}(t) $. 
Hence by \eqref{prf exc meas eq1} we obtain 
\begin{align}
\vn_{\rm BE}(e(A_m^{-1}(t+\cdot)) \in \Gamma ) 
=& \vn_{\rm BE}(e^+((A_m^+)^{-1}(\cdot)) \in \Gamma ) 
\\
=& \int_{(0,\infty )} \frac{1}{x} \vP^0_{\rm 3B}(A_m^{-1}(t) \in \d x) 
\vQ_{\rm BM}^x(e(A_m^{-1}(\cdot)) \in \Gamma ) 
\end{align}
for any measurable set $ \Gamma $ of $ E $. 
Therefore we obtain the formula {\rm (ii)}$ ' $ from Proposition \ref{prop: time-change}. 

The formula {\rm (i)}$ ' $ is an immediate consequence of \eqref{Qm formula}, 
which is a special case of the formula {\rm (ii)}$ ' $. 

The maximum decomposition formula {\rm (iii)}$ ' $ is obvious 
from Proposition \ref{prop: time-change}. 
In fact, the maximum value is invariant under time change, so $ M(e_m)=M(e) $. 
\end{proof}

\section{Proof of the convergence theorem} \label{sec: prf continuity}

We remark the following two elementary facts without proofs. 
\begin{Lem} \label{lem: unif conv}
Let $ f_{\lambda} $ and $ f $ be non-decreasing functions 
on $ [a,b] $ with $ -\infty<a<b<\infty $. 
Suppose that $ f $ is continuous 
and that $ f_{\lambda}(x) \to f(x) $ for all $ x \in [a,b] $. 
Then $ f_{\lambda}(x) \to f(x) $ uniformly in $ x \in [a,b] $. 
\end{Lem}

\begin{Lem} \label{lem: inverse}
Let $ f_{\lambda} $ and $ f $ be non-decreasing functions 
on $ [a,b] $ with $ -\infty<a<b<\infty $. 
Suppose that $ f $ is strictly increasing 
and that $ f_{\lambda}(x) \to f(x) $ as $ n \to \infty $ uniformly in $ x \in [a,b] $. 
Then $ f_{\lambda}^{-1}(x) \to f^{-1}(x) $ uniformly in $ x \in (f(a),f(b)) $. 
\end{Lem}

Now we proceed to prove Theorem \ref{thm: exc conv}. 
\begin{proof}[Proof of Theorem \ref{thm: exc conv}]
\setcounter{Ci}{0} 
\Circ 
By the definition \eqref{am}, 
it suffices to show that $ A_{m_{\lambda}}(t) \to A_m(t) $ as $ \lambda \to \infty $ 
for any fixed $ t>0 $ and for $ \vn_{\rm BE} $-a.e.~path $ e $. 
Let $ \delta>0 $. 
Then 
\begin{align}
\left| A_{m_{\lambda}}(t) - A_m(t) \right| 
\le {\rm (I)} + {\rm (II)} 
\label{conv 1}
\end{align}
where 
\begin{align}
{\rm (I)} = \left| \int_{(\delta,M(e)]} \lt(t,x) \d m_{\lambda}(x) 
- \int_{(\delta,M(e)]} \lt(t,x) \d m(x) \right| 
\end{align}
and 
\begin{align}
{\rm (II)} = \int_{(0,\delta]} \lt(\zeta ,x) \cbra{ \d m_{\lambda}(x) + \d m(x) } . 
\end{align}

\Circ 
By the pointwise convergence condition \eqref{ptwise conv}, 
we see that 
$ \int_{(a,b]} f(x) \d m_{\lambda}(x) \to \int_{(a,b]} f(x) \d m(x) $ 
as $ \lambda \to \infty $ 
for all $ 0<a<b<\infty $ and each bounded continuous function $ f $ on $ [a,b] $. 
Therefore $ {\rm (I)} $ converges to 0 as $ \lambda \to \infty $. 
\bigskip

\Circ 
Now to complete the proof it suffices to show that 
\begin{align}
\lim_{\delta \to 0+} \limsup_{\lambda \to \infty } {\rm (II)} = 0 
\qquad \text{$ \vn_{\rm BE} $-a.e.} 
\label{conv 3}
\end{align}
By an argument like that used in the proof of Lemma \ref{lem: am finite}, 
we can reduce the convergence \eqref{conv 3} to 
\begin{align}
\lim_{\delta \to 0+} \limsup_{\lambda \to \infty } 
\int_{(0,\delta]} \lt(\infty ,x) \cbra{ \d m_{\lambda}(x) + \d m(x) } 
= 0 
\qquad \text{$ \vP_{\rm 3B}^0 $-a.e.} 
\label{conv 4}
\end{align}
\bigskip

\Circ 
By the Ray--Knight theorem, we know that 
the process $ ( \lt(\infty ,x): x \ge 0 ) $ under $ \vP_{\rm 3B}^0 $ 
obeys the law of the two-dimensional Bessel-squared process starting from the origin, 
which we denote by $ \{ (U(x):x \ge 0) , \vU \} $. 
Note that the transition kernel of this process is given by $ q^2_x(a,b) \d b $ where 
\begin{align}
q^2_x(a,b) = \frac{1}{2x} \exp \cbra{- \frac{a+b}{2x}} I_0 \cbra{ \frac{\sqrt{ab}}{x} } 
, \qquad a,b \ge 0. 
\label{BESQ kernel}
\end{align}

The convergence \eqref{conv 4} is equivalent to 
\begin{align}
\lim_{\delta \to 0+} \limsup_{\lambda \to \infty } 
\int_{(0,\delta]} U(x) \cbra{ \d m_{\lambda}(x) + \d m(x) } = 0 
\qquad \text{$ \vU $-a.s.} 
\label{conv 5}
\end{align}
By the law of the iterated logarithm, 
there exists a finite random variable $ C $ such that 
$ U(x) \le C x \log \log(1/x) $ for $ x \in (0,1) $, $ \vU $-a.s. 
Hence we obtain \eqref{conv 5} 
by the assumption of the $ \cM_L $-tightness condition and the condition $ m \in \cM_L $. 
This completes  the proof of \eqref{am conv}. 
\bigskip

\Circ 
Now we apply Lemma \ref{lem: inverse} to obtain 
\begin{align}
\lim_{\lambda \to \infty } \sup_{t \ge 0} 
\left| A_{m_{\lambda}}^{-1}(t) - A_m^{-1}(t) \right| = 0 
\qquad \text{$ \vn_{\rm BE} $-a.e.} 
\end{align}
Since $ \vn_{\rm BE} $-a.e.~path $ e $ is uniformly continuous, 
we obtain the desired convergence \eqref{unif conv ae}. 
\end{proof}

Finally, we prove Proposition \ref{prop: exc conv converse}. 
\begin{proof}[Proof of Proposition \ref{prop: exc conv converse}]
We imitate the proof of Lemma 2 of \cite{MR1071537}. 
Assume that the convergence \eqref{am conv} holds. 
Let $ \delta>0 $ be fixed for a while. 
By the maximum decomposition formula \eqref{BE Williams}, we have 
\begin{align}
\lim_{\lambda \to \infty } \int_{(0,\delta]} \lt(\tau_{\delta},x) \d m_{\lambda}(x) 
= \int_{(0,\delta]} \lt(\tau_{\delta},x) \d m(x) 
, \ \forall \delta>0 \qquad \text{$ \vP^0_{\rm 3B} $-a.s.} 
\label{conv conv 1}
\end{align}

By the Ray--Knight theorem, we know that 
the process $ ( \lt(\tau_{\delta} ,x): x \in [0,\delta] ) $ under $ \vP_{\rm 3B}^0 $ 
obeys the law of the 2-dimensional Bessel-squared process starting from the origin 
pinned at the origin when $ x=\delta $, 
which we denote by $ \{ (U(x):x \in [0,\delta]) , \vU_{\delta} \} $. 
Note that 
\begin{align}
\vU_{\delta} \cbra{ U(x) \in \d a } 
=& \frac{q^2_x(0,a) q^2_{\delta-x}(a,0)}{q^2_{\delta}(0,0)} \d a 
\\
=& \frac{\delta}{2x(\delta-x)} \exp \cbra{- \frac{\delta a}{2x(\delta-x)} } \d a , 
\label{pinned BESQ}
\end{align}
where $ q^2_x(a,b) $ is given in \eqref{BESQ kernel}. 
Then we have 
\begin{align}
\vU_{\delta} \ebra{ U(x)/x } = \frac{2(\delta-x)}{\delta} \le 2 
, \qquad x \in (0,\delta] 
\label{conv conv 2}
\end{align}
and 
\begin{align}
\vU_{\delta} \cbra{ U(x)/x \le u } 
=& 1 - \exp \cbra{ -\frac{\delta u}{2(\delta -x)} } 
, \qquad x \in (0,\delta], \ u \in [0,\infty ) 
\\
\le& 1-\e^{-u} 
, \qquad x \in (0,\delta/2], \ u \in [0,\infty ) . 
\label{conv conv 3}
\end{align}
By the convergence \eqref{conv conv 1}, 
we can take $ \lambda_{\delta}>0 $ so large that 
$ \vU_{\delta}(B_{\delta}) \ge 1/2 $ where 
\begin{align}
B_{\delta} =& 
\kbra{ 
\lim_{\lambda \to \infty } \int_{(0,\delta]} U(x) \d m_{\lambda}(x) 
= \int_{(0,\delta]} U(x) \d m(x) } 
\\
& \cap 
\kbra{ 
\sup_{\lambda>\lambda_{\delta}} \int_{(0,\delta]} U(x) \d m_{\lambda}(x) \le L_{\delta} } 
\end{align}
and 
\begin{align}
L_{\delta} = \int_{(0,\delta]} U(x) \d m(x) + 1 . 
\end{align}
In fact, 
\begin{align}
\lim_{\lambda \to \infty } \vU_{\delta} \cbra{ \sup_{\lambda'>\lambda} 
\left| \int_{(0,\delta]} U(x) \d m_{\lambda'}(x) 
- \int_{(0,\delta]} U(x) \d m(x) \right| > 1 } = 0 . 
\end{align}
By \eqref{conv conv 2}, we have $ \vU_{\delta}[L_{\delta}] < \infty $. 
We can now apply the dominated convergence theorem to obtain 
\begin{align}
\lim_{\lambda \to \infty } \vU_{\delta} \ebra{ 1_B \int_{(0,\delta]} U(x) \d m_{\lambda}(x) } 
=& \vU_{\delta} \ebra{ 1_B \int_{(0,\delta]} U(x) \d m(x) } 
\\
\le& \int_{(0,\delta]} \vU_{\delta} \ebra{ U(x)/x } x \d m(x) 
\\
\le& 2 \int_{(0,\delta]} x \d m(x) . 
\label{conv conv 4}
\end{align}
On the other hand, we have 
\begin{align}
\vU_{\delta} \ebra{ 1_B \int_{(0,\delta]} U(x) \d m_{\lambda}(x) } 
=& \int_{(0,\delta]} \vU_{\delta} \ebra{ 1_B U(x)/x } x \d m_{\lambda}(x) 
\\
=& \int_{(0,\delta]} x \d m_{\lambda}(x) \int_0^{\infty } \vU_{\delta} \cbra{ \Big. B \cap \kbra{ U(x)/x > u } } \d u 
\\
\ge& \int_{(0,\delta]} x \d m_{\lambda}(x) \int_0^{\infty } \ebra{ \Big. \vU_{\delta}(B) - \vU_{\delta} \kbra{ U(x)/x \le u } }^+ \d u 
\\
\ge& C \int_{(0,\delta/2]} x \d m_{\lambda}(x) 
\label{conv conv 5}
\end{align}
where $ C = \int_0^{\infty } \ebra{ \e^{-u} - 1/2 }^+ \d u > 0 $. 
(Here we used \eqref{conv conv 3}.)
Therefore we obtain 
\begin{align}
\limsup_{\lambda \to \infty } \int_{(0,\delta/2]} x \d m_{\lambda}(x) 
\le \frac{2}{C} \int_{(0,\delta]} x \d m(x) . 
\label{conv conv 6}
\end{align}
If we let $ \delta $ tend to $ 0+ $, 
then the RHS of \eqref{conv conv 6} vanishes. 
\end{proof}

\section{Proof of the limit theorem} \label{sec: prf lim thm}

Firstly, we prove Lemma \ref{lem: lim thm}. 
\begin{proof}[Proof of Lemma \ref{lem: lim thm}]
\setcounter{Ci}{0} 
\Circ 
In the case of {\bf (A1)}, 
the inequality \eqref{global bound} is easily justified 
since $ 0 < m_{\lambda}(x) \le m_{\lambda}(1) $ for $ x \in (0,1] $ and for $ \lambda>0 $ 
and $ m_{\lambda}(1) \to 1 $ as $ \lambda \to \infty $. 

\Circ 
In the case of {\bf (A2)}, 
we may take $ m(0):=m(0+) $ so that $ m(x) $ is locally bounded on $ [0,\infty ) $. 
Hence we may apply Theorem 3.8.6 (b) of \cite{MR1015093}, pp.172,
and then we obtain the inequality \eqref{global bound}. 

\Circ 
In the case of {\bf (A3)}, 
using the two conditions \eqref{A3 prime at 0} and \eqref{A3 prime at infty}, 
we may take a constant $ C_1 $ 
and a function $ \tilde{K}(x) $ defined on $ [0,\infty ) $ 
such that the following hold: 
\\
{\rm (i)} 
$ m(\infty ) - m(x) \le C_1 x^{\frac{1}{\alpha }-1} \tilde{K}(x) $ 
for all $ x>0 $. 
\\
{\rm (ii)} 
$ \tilde{K}(x) $ is bounded away from 0 and $ \infty $ on any compact subset of $ [0,\infty ) $. 
\\
{\rm (iii)} 
$ \tilde{K}(x)/K(x) \to 1 $ as $ x \to \infty $. 
(Then $ \tilde{K}(x) $ is necessarily slowly varying as $ x\to\infty $.)

We may apply Theorem 1.5.6 (ii) of \cite{MR1015093}, pp.25, 
to the function $ \tilde{K}(x) $, 
and see that there exist a constant $ C $ such that 
\begin{align}
\tilde{K}(\lambda x)/\tilde{K}(\lambda) \le C_2 x^{- \frac{1}{2 \alpha }-1} , 
\qquad x \in (0,1] , \ \lambda>0 . 
\end{align}
Now we obtain 
\begin{align}
|m_{\lambda}(x)| 
= \frac{m(\infty ) - m(\lambda x)}{\lambda^{\frac{1}{\alpha }-1} K(\lambda)} 
\le C x^{\frac{1}{2 \alpha }-1} 
\qquad x \in (0,1] , \ \lambda>0 
\end{align}
for some constant $ C $. 
Hence we obtain the inequality \eqref{global bound} with $ \eps = \frac{1}{2 \alpha } $. 

\Circ 
Since $ \int_0^{1/2} x^{\eps-1} \log \log (1/x) \d x < \infty $, 
we can apply the dominated convergence theorem to obtain 
\begin{align}
\lim_{\lambda \to \infty } \int_0^{\delta} |m_{\lambda}(x)| \log \log(1/x) \d x 
= \int_0^{\delta} |m^{(\alpha )}(x)| \log \log(1/x) \d x , 
\qquad \delta < 1/2 . 
\label{dominate limit}
\end{align}
Therefore we conclude that the $ \cM_L $-tightness condition \eqref{cMlog2 tight} holds 
by \eqref{dominate} and \eqref{dominate limit}. 
Since it is obvious that $ m_{\lambda}(x) \to m^{(\alpha )}(x) $ pointwise, 
we conclude that $ m_{\lambda} \to m $ in $ \cM_L $. 
\end{proof}

Secondly, we prove Lemma \ref{lem: lim thm2}. 
\begin{proof}[Proof of Lemma \ref{lem: lim thm2}]
By \eqref{Qm formula}, we have 
\begin{align}
\vQ_m^x \ebra{ F \cbra{ \lambda^{-1} e(v(\lambda) t ) : t \ge 0 } }
= x 
\vn_m \ebra{ F \cbra{ \lambda^{-1} e(v(\lambda) t + \tau_x ) : t \ge 0 } } . 
\label{lem: lim thm2 eq2}
\end{align}
Since $ \tau_x(e_m) = A_m(\tau_x) $ and $ \zeta(e_m) = A_m(\zeta) $, the RHS is equal to 
\begin{align}
x 
\vn_{\rm BE} \ebra{ 
F \cbra{ \lambda^{-1} e( A_m^{-1}( v(\lambda) t + A_m(\tau_x) ) ) : t \ge 0 } } . 
\label{lem: lim thm2 eq3}
\end{align}
Now consider $ \tilde{e}(t) = \lambda e(t/\lambda^2) $. 
By the scaling property of Bessel processes, 
we have $ \vn_{\rm BE}(\tilde{e} \in \cdot) 
= \lambda \vn_{\rm BE}(e \in \cdot) $. 
Then \eqref{lem: lim thm2 eq3} can be rewritten as 
\begin{align}
\frac{x}{\lambda} 
\vn_{\rm BE} \ebra{ F \cbra{ \lambda^{-1} \tilde{e} \cbra{ \tilde{A_m}^{-1}( v(\lambda) t + \tilde{A_m}(\tilde{\tau_x}) ) } : t \ge 0 } } . 
\label{lem: lim thm2 eq4}
\end{align}
Here note that $ \tilde{\zeta}=\zeta(\tilde{e})=\lambda^2 \zeta $, 
that $ \tilde{\tau_x} = \tau_x(\tilde{e}) = \lambda^2 \tau_{x/\lambda} $, 
and that $ \tilde{A_m}(t) = v(\lambda) A_{m_{\lambda}}(t/\lambda^2) $. 
Thus \eqref{lem: lim thm2 eq4} leads to 
\begin{align}
\frac{x}{\lambda} 
\vn_{\rm BE} \ebra{ F \cbra{ e \cbra{ A_m^{-1}( t + A_m(\tau_{x/\lambda}) ) } : t \ge 0 } } . 
\label{lem: lim thm2 eq5}
\end{align}
If we replace $ F(e) $ by $ F(e) 1_{\{ \zeta(e)>1 \}} $, then 
we obtain \eqref{lem: lim thm2 eq1} 
from \eqref{lem: lim thm2 eq2} and \eqref{lem: lim thm2 eq5}. 
\end{proof}

Thirdly, we prove Lemma \ref{lem: lim thm3}. 
\begin{proof}[Proof of Lemma \ref{lem: lim thm3}]
\setcounter{Ci}{0} 
\Circ 
By Theorem \ref{thm: exc conv}, 
we have $ F(e_{m_{\lambda}}) \to F(e_m) $ as $ \lambda \to \infty $. 
By Theorem \ref{thm: exc conv} again, 
we have $ A_{m_{\lambda}}(\zeta) - A_{m_{\lambda}}(\tau_{x/\lambda}) \to A_m(\zeta) $ 
$ \vn_{\rm BE} $-a.e. 
Since $ \vn_{\rm BE}(A_m(\zeta)=1)=0 $, we see that 
the indicator function of 
$ \{ A_{m_{\lambda}}(\zeta) - A_{m_{\lambda}}(\tau_{x/\lambda}) > 1 \} $ 
converges to that of 
$ \{ A_m(\zeta) > 1 \} $ 
$ \vn_{\rm BE} $-a.e. 
Hence, for every $ \delta>0 $, we obtain 
\begin{align}
\begin{split}
 \lim_{\lambda \to \infty } \vn_{\rm BE}& \ebra{ F(e_{m_{\lambda}}) ; A_{m_{\lambda}}(\zeta) - A_{m_{\lambda}}(\tau_{x/\lambda}) > 1 , \ M>\delta } 
\\
&= \vn_{\rm BE} \ebra{ F(e_m) ; A_m(\zeta) > 1 , \ M>\delta } 
\end{split}
\label{lem: lim thm3 eq2}
\end{align}
by the dominated convergence theorem. 

\Circ 
Now it is sufficient to show that 
\begin{align}
\lim_{\delta \to 0+} \sup_{\lambda>\lambda_0} 
\vn_{\rm BE} \ebra{ F(e_{m_{\lambda}}) ; A_{m_{\lambda}}(\zeta) - A_{m_{\lambda}}(\tau_{x/\lambda}) > 1 , \ M \le \delta } = 0 . 
\end{align}
Noting that $ F $ is bounded 
and that $ A_{m_{\lambda}}(\zeta) - A_{m_{\lambda}}(\tau_{x/\lambda}) 
\le A_{m_{\lambda}}(\zeta) $, 
we have only to show that 
\begin{align}
\lim_{\delta \to 0+} \sup_{\lambda>\lambda_0} 
\vn_{\rm BE} \cbra{ A_{m_{\lambda}}(\zeta) > 1 , \ M \le \delta } = 0 . 
\end{align}
Using the terminology in the proof of Lemma \ref{lem: am finite}, 
we have $ \{ A_{m_{\lambda}}(\zeta) > 1 \} \subset \{ A_{m_{\lambda}}(\tau_a)>1/2 \} \cup \{ A_{m_{\lambda}}^{\vee}(\tau_a^{\vee})>1/2 \} $, 
and hence we have 
\begin{align}
\vn_{\rm BE} \cbra{ A_{m_{\lambda}}(\zeta) > 1 , \ M \le \delta } 
\le 2 \int_0^{\delta} \vP_{\rm 3B}^0 (A_{m_{\lambda}}(\tau_a)>1/2) \frac{\d a}{a^2} 
\label{lem: lim thm3 eq3}
\end{align}
by the maximum decomposition formula \eqref{BE Williams}. 
Noting that $ A_{m_{\lambda}}(\tau_a) \le \int_{(0,a]} \lt(\infty ,x) \d m_{\lambda}(x) $ 
and that $ (\lt(\infty,x):x \ge 0) $ under $ \vP_{\rm 3B}^0 $ 
has the same law as the two-dimensional  Bessel-squared  process, 
the RHS of \eqref{lem: lim thm3 eq3} is dominated by 
\begin{align}
2 \int_0^{\delta} \vP \cbra{ \int_{(0,a]} B(x)^2 \d m_{\lambda} (x) >1/4 } \frac{\d a}{a^2} 
\label{lem: lim thm3 eq4}
\end{align}
where $ \{ (B(x):x \ge 0) , \vP \} $ is the law of the standard Brownian motion. 
Hence it suffices to find a function $ f(a) $ defined near $ a=0 $ 
with $ \int_{0+} f(a) \frac{\d a}{a^2} < \infty $ 
such that 
$ \sup_{\lambda>\lambda_0} 
\vP ( \int_{(0,a]} B(x)^2 \d m_{\lambda} (x) >1/4 ) \le f(a) $ near $ a=0 $. 

\Circ 
Now we proceed to find such a function $ f(a) $. 
Set $ m^a_{\lambda} (x) = a^{1-\eps} m_{\lambda} (ax) $. 
Then 
\begin{align}
|m^a_{\lambda} (x)| \le C x^{\eps -1} 
\qquad \text{for all $ x,a \in (0,1] $ and $ \lambda>\lambda_0 $} 
\label{lem: lim thm3 eq5}
\end{align}
by the assumption \eqref{global bound}. 
By the scaling property of Brownian motion, 
the random variable $ \int_{(0,a]} B(x)^2 \d m_{\lambda} (x) $ 
has the same law as 
$ a^{\eps} \int_{(0,1]} B(x)^2 \d m^a_{\lambda} (x) $ 
under $ \vP $. 
By  Chebyshev's inequality, we have 
\begin{align}
\vP \cbra{ a^{\eps} \int_{(0,1]} B(x)^2 \d m^a_{\lambda} (x) >1/4 } 
\le \e^{ - \rho a^{-\eps}/4 } 
\vP \ebra{ \exp \kbra{ \rho \int_{(0,1]} B(x)^2 \d m^a_{\lambda} (x) } } 
\end{align}
for all $ a \in (0,1] $ and for all $ \rho>0 $. 
Hence it suffices to find $ \rho>0 $ such that 
\begin{align}
C':= 
\sup_{\lambda>\lambda_0} \sup_{a \in (0,1]} 
\vP \ebra{ \exp \kbra{ \rho \int_{(0,1]} B(x)^2 \d m^a_{\lambda} (x) } } 
< \infty . 
\end{align}
Indeed, we can choose $ f(a)=C' \e^{ - \rho a^{-\eps}/4 } $. 

\Circ 
By It\^o's formula, we see that 
\begin{align}
\int_{(0,1]} B(x)^2 \d m^a_{\lambda} (x) = G_{\lambda}^a + H_{\lambda}^a + I_{\lambda}^a 
\label{lem: lim thm3 eq6}
\end{align}
where $ G_{\lambda}^a = B(1)^2 m^a_{\lambda}(1) $, 
$ I_{\lambda}^a = - \int_0^1 m^a_{\lambda} (x) B(x)^2 \d x $ 
and 
\begin{align}
H_{\lambda}^a = - 2 \int_0^1 m^a_{\lambda} (x) B(x) \d B(x) . 
\end{align}
Here we remark that 
$ \lim_{x \to 0+} B(x)^2 |m^a_{\lambda}(x)| = 0 $ $ \vP $-a.s. by \eqref{lem: lim thm3 eq5}. 
Note that $ G_{\lambda}^a \le C B(1)^2 $ 
and that $ I_{\lambda}^a \le C \int_0^1 x^{\eps-1} B(x)^2 \d x $. 
We also note that $ M(x) = -2 \int_0^x m^a_{\lambda} (y) B(y) \d B(y) $ is a martingale 
such that 
\begin{align}
\langle M \rangle (1) = 4 \int_0^1 m^a_{\lambda} (x)^2 B(x)^2 \d x 
\le 4 C^2 \int_0^1 x^{2\eps -2} B(x)^2 \d x . 
\end{align}
Since $ B(1)^2 $, $ \int_0^1 x^{\eps-1} B(x)^2 \d x $ 
and $ \int_0^1 x^{2\eps -2} B(x)^2 \d x $ 
are all quadratic Wiener functionals, they are exponentially integrable 
if we take the exponent sufficiently small. 
\end{proof}

\begin{proof}[Proof of Theorem \ref{thm: lim thm}]
Since the law of the lifetime of the Bessel meander is absolutely continuous, 
we have $ \vn_{\rm BE}(A_{m^{(\alpha )}}(\zeta) = 1) = 0 $. 
Now we apply Lemmas \ref{lem: lim thm}, \ref{lem: lim thm2} and \ref{lem: lim thm3},
putting $ v(\lambda) = \lambda^{\frac{1}{\alpha }} K(\lambda) $, to obtain the desired conclusion. 
\end{proof}

\end{document}